\documentclass[a4paper,oneside,11pt]{article}
\usepackage{amsmath,amsthm,amsfonts,amssymb}

\theoremstyle{plain}
\newtheorem{Thm}{Theorem}
\newtheorem{Lem}[Thm]{Lemma}

\newtheorem{Cor}[Thm]{Corollary}

\theoremstyle{definition}
\newtheorem{Defn}[Thm]{Definition}

\newenvironment{pf}{ \begin{proof} }{ \end{proof} }

\newcommand{\R}{\mathbb{R}}
\newcommand{\C}{\mathbb{C}}

\DeclareMathOperator{\sym}{\mathrm{Sym}}

\begin{document}
\author{Tim Perutz}
\title{A remark on K\"ahler forms on symmetric products of Riemann surfaces}
\date{January 31, 2005}


\maketitle

\begin{abstract}
Users of Heegaard Floer homology may be reassured to know that it can be made to conform exactly to the standard analytic pattern of Lagrangian Floer homology. This follows from the following remark, which we prove using an argument of J. Varouchas: the natural singular K\"ahler form $\sym^n(\omega)$ on the $n$th symmetric product of a K\"ahler curve $(\Sigma,\omega)$ admits a cohomologous smoothing to a K\"ahler form which equals $\sym^n(\omega)$ away from a chosen neighbourhood of the diagonal. 
\end{abstract}

In this note we consider branched coverings $\pi\colon X\to X' $ of complex manifolds---that is, holomorphic maps which are proper, surjective, and finite. The branch locus $B_\pi \subset X'$ of such a map is 
\[ B_\pi = \{ \pi(x): x\in X, \; \ker D_x(\pi) \neq 0\}.\] 
A ($C^\infty$) smooth K\"ahler form $\omega$ on $X$ can be pushed forward---in the sense of currents, that is, of 2-forms with $L^1$ coefficients---to a closed current $\pi_*(\omega)$ on $X'$ which is smooth on $X'\setminus B_\pi$. 

The following theorem is essentially due to Varouchas \cite{Var}.
\begin{Thm}\label{var}
Let $\pi\colon  X \to X'$ be a branched covering of complex manifolds,
and $\omega$ a K\"ahler form on $X$. Let $N$ be a
neighbourhood of the branch locus in $X'$. Then there exists a
K\"ahler form $\omega'$ on $X'$ such that
\begin{enumerate}
\item   $(\pi_*\omega-\omega')|_{X'\setminus \overline{N}}=0$, and
\item   $[\omega']=\pi_*[\omega]\in H^2(X';\R)$.
\end{enumerate}
\end{Thm}
The stated conclusion in \cite{Var} is simply that $X'$ admits a K\"ahler form. The purpose of this note is to explain this minor modification of Varouchas' argument, and to draw attention to the following example: 
\begin{Cor}
Let $\Sigma$ be a Riemann surface with K\"ahler form $\omega$. Let 
\[  \pi \colon \Sigma^{\times r}\to \sym^r(\Sigma) = \Sigma^{\times r}/S_r \]
be the projection map. Suppose that $N\subset \sym^r(\Sigma)$ is an open subset containing the (large) diagonal. Then there exists a K\"ahler form $\eta$ on $\sym^r(\Sigma)$ such that
\begin{enumerate}
\item
outside $\overline{N}$, $\eta$ is the smooth push-forward $\pi_*(\omega^{\times r})$ of the product form;
\item
$[\eta]=\pi_*[\omega^{\times r}]\in H^2(\sym^r(\Sigma);\R)$.
\end{enumerate}
\end{Cor}

As advertised in the abstract, this has a direct application to the Heegaard Floer theory of Ozsv\'ath-Szab\'o \cite{OS}. There one has tori $\mathbb{T}_\alpha$ and $\mathbb{T_\beta}$ in $\sym^g(\Sigma)$, disjoint from the diagonal and Lagrangian for product forms $\omega^{\times r}$; they remain Lagrangian for a globally-defined K\"ahler form $\eta$. The construction of Ozsv\'ath-Szab\'o differs from the standard analytic framework of Lagrangian Floer homology only in the handling of energy bounds for the holomorphic disks; if one uses the form $\eta$, the energy has its usual cohomological interpretation, and so is constant on each homotopy class.

\vspace{3mm}

{\bf Acknowledgements.} \emph{This note was drafted in 2003, with different applications in mind; I apologise for the delay in releasing it. Thanks to Simon Donaldson and Sebastien Boucksom for useful discussions at the time, and to Alexandru Oancea for helpful comments.}

\vspace{3mm}

\begin{Defn}
Let $X$ be a complex manifold. A {\bf K\"ahler cocycle} on $X$ is
a collection $(U_i,\varphi_i)_{i\in I}$, where $(U_i)_{i\in
I}$ is an open cover of $X$, and $\varphi_i\colon U_i \to \R$ is a
a function, such that for all $i,j \in I$,
\begin{enumerate}
\item $\varphi_i$ is strictly plurisubharmonic on $U_i$; and
\item $\varphi_i-\varphi_j$ is pluriharmonic on $U_i \cap U_j$.
\end{enumerate}
One ascribes to the cocycle a property (continuity, smoothness,
etc.) possessed by all the $\varphi_i$. K\"ahler cocycles are, by
definition, upper semicontinuous.
\end{Defn}
Condition (1) means that the $2$-current $dd^c \varphi_i$ is
strictly positive on $U_i$; (2) means that these currents agree on
overlaps, and are therefore restrictions of a $2$-current $\omega$
on $X$ (closed and strictly positive). If the cocycle is
$C^\infty$ then $\omega$ will be a K\"ahler form.

Varouchas' \emph{lemme principal} is the following. The proof uses the ``regularised maximum" technique of Richberg and Demailly.

\begin{Lem}
Let $U,V,W,\Omega$ be open subsets of $\C^n$ with
\[U\Subset V\Subset W, \quad \Omega \subset W. \]
Let $\phi\colon W \to \R$ be continuous, strictly plurisubharmonic, and smooth on $\Omega$. 
Then there exists a function $\psi \colon W\to \R$, again continuous and strictly plurisubharmonic, equal to $\phi$ on $W \setminus \overline{V}$ and smooth on $U\cup \Omega$.
\end{Lem}
One then passes from local to global by the following argument, which I give in detail since Varouchas' stated conclusion is weaker here. 
\begin{Lem}
Let $(U_i,\varphi_i)_{i\in I}$ be a continuous K\"ahler cocycle on the complex manifold $X$. Suppose that $X=X_1\cup X_2$, with $X_1$ and $X_2$ open,
and that the functions $\varphi_i|_{U_i \cap X_1}$ are smooth. Then
there exists a continuous function
\[\chi \colon X \to \R, \quad \mathrm{Supp}(\chi) \subset X_2,\]
and a locally finite refinement
\[ V_j \subset U_{i(j)} \quad (j \in J) \]
so that the family
\[   (V_j, \varphi_{i(j)}|_{V_j}+\chi|_{V_j})_{j\in J}     \]
is a smooth K\"ahler cocycle.
\end{Lem}
\begin{pf}
Refine the cover $(U_i)_{i\in I}$ to a countable, locally finite
cover $(V_i)_{i\in I_1 \coprod I_2}$ with the property that
\[ i \in I_\alpha \;  \Rightarrow  \;V_i \subset X_\alpha, \quad \alpha=1,2. \]
For definiteness let us suppose both $I_1$ and $I_2$ are infinite;
say $I_\alpha = \mathbb{N}\times \{\alpha\}.$ Find open subsets
\[V_i'' \Subset V_i' \Subset V_i \]
such that $(V_i'')$ still covers $X$, and set
\[  A_1=\emptyset, \quad A_n = V''_{(1,2)}\cup\dots\cup V''_{(n-1,2)}. \]
So the sets $A_n$ exhaust $X_2\setminus X_1$. Let
$(V_i,\psi^1_i)_{i\in I_1\cup I_2}$ be the K\"ahler cocycle
induced from $(U_i,\varphi_i)_{i\in I}$ by the
refinement.

\emph{Claim:} there are K\"ahler cocycles $(V_i,\psi^n_i)$, where $n=1,2,\dots$ indexes the elements of
$I_2$, such that the following hold for all $i \in I_1 \cup I_2$ and all $n>1$:
\begin{enumerate}
\item   $\psi^n_i$ is smooth on the set $V_i \cap (X_1\cup A_n).$
\item   There is a continuous function $\chi_n\colon X \to \R$, with $\mathrm{Supp}(\chi_n)\subset
        V'_{(n-1,2)}$, such that $\psi^n_i=  \psi^{n-1}_i +\chi_n$.
\end{enumerate}
We prove the claim by induction on $n$. Apply the previous lemma to
\[ (U,V,W,\Omega) = (V''_{(n-1,2)},V'_{(n-1,2)},V_{(n-1,2)}, V_{(n-1,2)}\cap(X_1\cup A_{n-1})) \]
and to the function $\psi^{n-1}_{(n-1,2)}$, obtaining a new function $\psi^{n}_{(n-1,2)} $; let $\chi_n = \psi^{n}_{(n-1,2)} - \psi^{n-1}_{(n-1,2)}$, extended by zero to all of $X$, and use (2) to define the new cocycle. We have to verify (1), i.e. to prove smoothness of $\psi^n_i$ at each $x\in V_i\cap(X_1 \cup A_{n})$. If $x\not \in V'_{(n-1,2)}$ then $\chi_n(x)=0$, but $\psi_i^{n-1}$ was already smooth. If $x\in V_{(n-1,2)}'$ then, near $x$, $\psi^n_i = (\psi^n_i -\psi^n_{(n-1,2)})+ \psi^n_{(n-1,2)} = (\psi^1_i - \psi^1_{(n-1,2)}) + \psi^n_{(n-1,2)}$, which is the sum of a pluriharmonic function and a smooth plurisubharmonic one. But a pluriharmonic function is smooth. By a similar argument, $\psi^n_i$ is strictly plurisubharmonic.

Now define a function $\chi\colon X \to \R$ by the locally finite sum
\[ \chi(x) = \sum_{n \geq 1}{\chi_n(x)}.\]
Then $ \psi^\infty_i(x)  :=  \psi^1_i(x) +\chi(x)$
defines a K\"ahler cocycle. It is
smooth, since on $X \setminus {\bigcup{V'_{(n,2)}}} \subset X_1$, the original
cocycle was smooth and has not been modified, while
\[ V_{(n,2)} \subset X_1 \cup \bigcup{A_k},  \]
so smoothness on $V_{(n,2)}$ is guaranteed by (1). Hence $\chi$ has the
required properties.
\end{pf}
\begin{pf}[Proof of Theorem \ref{var}]
Each fibre $\pi^{-1}(x')$, being finite, has a neighbourhood which is a disjoint union of open balls. Hence, using the $dd^c$-lemma, one can find a smooth K\"ahler cocycle $(U_i,\varphi_i)$ on $X$ such that each $U_i$ contains a fibre of $\pi$, with $\omega|_{U_i}=dd^c \varphi_i$. 
One can then find a locally finite cover $(U'_i)$ of $X'$ such that
$U_i \supset \pi^{-1}(U'_i)$.

A general property of branched covers is that the push-forward $\pi_*f $ of a continuous function $f\colon X\to \R$ is again continuous (it is given by $\pi_*f (x')=\sum_{x\in\pi^{-1}(x')}f(x)$, where the points $x$ are taken with multiplicities). The family $(U'_i,\pi_*\varphi_i)$ on $X'$ is thus a continuous K\"ahler cocycle: plurisubharmonicity is clear away from $B_\pi$, hence
everywhere by density; similarly for pluriharmonicity on overlaps.

Now let $N' \Subset N$ be a smaller open neighbourhood of $B_\pi$.
Apply the global smoothing lemma to $(U'_i,\pi_*\varphi_i)$ on
$X'$, taking $X_1=X'\setminus \overline{N'}$ and $X_2=N$. The
output is a function $\chi\colon X' \to \R$ (as well as a
refinement of $(U_i')$, which we omit from the notation) such that
\[  \omega_{X'}:= dd^c (\pi_*\varphi_i +\chi)= \pi_*dd^c\varphi_i+dd^c\chi \]
is a well-defined $2$-form with the right properties. Notice that, since $\chi$ is continuous, $dd^c \chi$ represents the zero cohomology class.
\end{pf}

\flushleft{
Department of Mathematics, 

Imperial College London, 

London SW7 2AZ, UK. 

\vspace{2mm}

tim.perutz@imperial.ac.uk}

\end{document}